\documentclass[letterpaper, 10 pt, conference]{ieeeconf}

\usepackage{graphicx}
\usepackage{amsmath, amssymb, mathtools}
\usepackage{xcolor}
\usepackage{clevethm}
\usepackage{caption}
\usepackage{subfig}
\usepackage{tikz}
\usepackage{circuitikz}
\usepackage{pgfplots}
\usepackage{siunitx}
\usepackage{hyperref}

\IEEEoverridecommandlockouts
\DeclareMathOperator{\graph}{gph}
\DeclareMathOperator{\realpart}{Re}
\DeclareMathOperator{\impart}{Im}
\DeclareMathOperator{\dom}{dom}
\DeclareMathOperator{\ran}{ran}
\DeclareMathOperator{\id}{id}
\DeclareMathOperator{\fix}{Fix}
\DeclareMathOperator{\zer}{Zer}
\DeclareMathOperator{\rank}{rank}
\DeclareMathOperator{\adj}{adj}
\DeclareMathOperator{\arcsinh}{arcsinh}

\newcommand{\KM}{Krasnosel'ski\v{\i}--Mann}

\newcommand{\proj}{\Pi}

\newcommand{\cH}{\mathcal{H}}
\newcommand{\cS}{\mathcal{S}}
\newcommand{\cP}{\mathcal{P}}

\newcommand{\cM}{\mathcal{M}}
\newcommand{\cA}{\mathcal{A}}
\newcommand{\cG}{\mathcal{G}}
\newcommand{\bR}{\mathbb{R}}
\newcommand{\bC}{\mathbb{C}}

\newenvironment{iproof}{
  \proof}{\endproof}

\newcommand{\sing}{F}

\newtheorem{theorem}{Theorem}[section]
\newtheorem{lemma}{Lemma}[section]
\newtheorem{proposition}{Proposition}[section]
\newtheorem{corollary}{Corollary}[section]
\newtheorem{definition}{Definition}[section]
\theoremstyle{remark}
\newtheorem{remark}{Remark}[section]
\newtheorem{example}{Example}[section]

\crefname{definition}{Definition}{Definitions}
\Crefname{definition}{Definition}{Definitions}
\crefname{theorem}{Theorem}{Theorems}
\Crefname{theorem}{Theorem}{Theorems}
\crefname{lemma}{Lemma}{Lemmas}
\Crefname{lemma}{Lemma}{Lemmas}
\crefname{proposition}{Proposition}{Propositions}
\Crefname{proposition}{Proposition}{Propositions}
\crefname{example}{Example}{Examples}
\Crefname{example}{Example}{Examples}
\crefname{corollary}{Corollary}{Corollaries}
\Crefname{corollary}{Corollary}{Corollaries}
\crefname{figure}{Fig.}{Figs.}
\Crefname{figure}{Fig.}{Figs.}

\newlist{propenum}{enumerate}{1} 
\setlist[propenum]{label=(\roman*), ref=\theproposition(\roman*), font=\rm}
\crefalias{propenumi}{proposition}

\hypersetup{
hypertexnames = false,
plainpages= false,
colorlinks= true,
citecolor = ForestGreen,
linkcolor = NavyBlue,
}

\begin{document}

\title{
   Scaled relative graphs for pairs of operators beyond classical monotonicity
        \thanks{
            This work was supported by the Research Foundation Flanders (FWO) PhD grant 11A8T26N and research projects G081222N, G033822N, and G0A0920N; Research Council KUL grant C14/24/103.
        }
}

\author{%
	\texorpdfstring{%
		Jan Quan\thanks{
			KU Leuven, Department of Electrical Engineering ESAT-STADIUS -- %
			Kasteelpark Arenberg 10, box 2446, 3001 Leuven, Belgium
			{\sf\scriptsize
					\href{mailto:jan.quan@kuleuven.be}{jan.quan@kuleuven.be}%
			}%
		}
		\and
        Alexander Bodard
		\and
		Konstantinos Oikonomidis
		\and
		Panagiotis Patrinos
	}{%
		Jan Quan, Alexander Bodard, Konstantinos Oikonomidis and Panagiotis Patrinos
	}%
}

\Crefname{figure}{Fig.}{Figs.}

\setlength{\textfloatsep}{4pt plus 2pt minus 2pt}      

\setlength{\floatsep}{3pt plus 2pt minus 2pt}      
\setlength{\intextsep}{6pt plus 2pt minus 2pt}

\maketitle
\begin{abstract}
We introduce a generalization of the scaled relative graph (SRG) to pairs of operators, enabling the visualization of their relative incremental properties. This novel SRG framework provides the geometric counterpart for the study of nonlinear resolvents based on paired monotonicity conditions. We demonstrate that these conditions apply to linear operators composed with monotone mappings, a class that notably includes NPN transistors, allowing us to compute the response of multivalued, nonsmooth, and highly nonmonotone electrical circuits.
\end{abstract}

\section{Introduction}
Understanding the input-output behavior of systems and their underlying operators is a central problem in the study of dynamical systems and control algorithms. In this context, classical monotone operator theory provides a unifying mathematical framework for modeling feedback interconnections, optimization dynamics, and equilibrium systems. Moreover, monotonicity is a fundamental property for ensuring stability and convergence in many settings, from proximal algorithms to feedback systems governed by maximal monotone mappings. However, many relevant control and learning systems are described by nonmonotone operators and thus require a different approach.

A new concept, \emph{pair monotonicity}, which is a specialization of the $T$-monotonicity concept from \cite[Def.\,3]{borwein2011characterization}, has been introduced in order to better analyze the differential inclusions of sweeping processes in \cite{adly2024state}. Notably, a set-valued variant of this monotonicity property was also used in \cite{le2025solving} to characterize the firm nonexpansiveness of generalized resolvent operators. Departing from standard monotonicity that links the inputs and outputs of an operator, pair monotonicity describes the incremental properties of the output of one operator compared to the output of another, thus making it useful for characterizing highly nonmonotone systems. Pair monotonicity provides a powerful way to incorporate prior knowledge by choosing the second operator judiciously. As this approach is purely algebraic, this choice may be quite difficult in practice.

Rather than relying on algebra, graphical tools are invaluable for building intuition. To this end, the scaled relative graph (SRG) \cite{ryu2022scaled} has emerged as a powerful framework by mapping the action of operators onto the (extended) complex plane. Serving as a nonlinear generalization of the classical Nyquist diagram, this tool has since seen applications in various systems and control contexts, such as graphical system analysis \cite{chaffey2023graphical,baron2025stability,de2025exploiting,krebbekx2025graphical} and reset control systems~\cite{van2024scaled}, and has proven useful in the convergence analysis of various algorithms \cite{lee2025convergence,ryu2022scaled}.

Originally, the SRG was used for the analysis of operator properties, where the high-level approach is shown in \Cref{fig:meta_proc}. Of particular interest is showing that operators are firmly nonexpansive or contractive, upon which convergence of the associated fixed-point iteration can be shown using the \KM\ \cite[Cor.\,5.17]{bauschke2017correction} and Banach fixed-point theorem \cite[Thm.\,1.50]{bauschke2017correction}, respectively. With the emergence of the new pair monotonicity, it is natural to ask how the SRG can be extended to handle these novel properties. This is the main subject of study in this paper.

\begin{figure}[hbpt]
    \centering
    \begin{tikzpicture}
        \draw (0, 0) rectangle (2, 1);
        \node[anchor=center] at (1, 0.5) {$A\in\cA$};
        \draw (4.5, 0) rectangle (6.5, 1);
        \node[anchor=center] at (5.5, 0.5) {$B\in\mathcal{B}$};
        \draw (0, -1) rectangle (2, -2);
        \node[anchor=center] at (1, -1.5) {$\cG(A, \id)$};
        \draw (4.5, -1) rectangle (6.5, -2);
        \node[anchor=center] at (5.5, -1.5) {$\cG(B, \id)$};

        \draw[-stealth, thick] (2, 0.5) -- node[midway, above] {Algebra} (4.5, 0.5);
        \draw[-stealth, thick] (2, -1.5) -- node[midway, below] {Geometry} (4.5, -1.5);
        \draw[-stealth, thick] (1, 0) --node[midway, left] {$\cG(\cdot, \id)$} (1, -1);
        \draw[-stealth, thick] (5.5, -1) --node[midway, right] {SRG-full} (5.5, 0);
    \end{tikzpicture}
    \caption{We can derive properties of $B$ from properties of $A$ by transforming the SRG of $A$ (denoted as $\mathcal{G}(A,\id)$) into the SRG of $B$, and then using SRG-fullness. This offers a geometric alternative to the classical algebraic approach.}
    \label{fig:meta_proc}
\end{figure}

Concretely, our contribution is threefold.
\begin{enumerate}[label=(\roman*)]
    \item We introduce a scaled relative graph for pairs of operators and establish important calculus rules. Furthermore, we extend the notions of SRG-fullness and semimonotonicity to this setting, thereby extending classical monotonicity.
    \item We apply this novel tool to provide purely geometric proofs for the core properties of two nonlinear resolvents that have been used to solve inclusion problems with pairs of monotone operators.
    \item We show the practical utility of this paired monotonicity condition by analyzing operators that can be written as the composition of a linear operator with a monotone mapping, a class that includes nonlinear transistor models. Leveraging this property, we are able to compute the response of a nonmonotone common-emitter amplifier circuit.
\end{enumerate}

\subsection{Notation}
In the following, $\cH$ denotes a real Hilbert space with inner product $\langle \cdot,\cdot\rangle$ and induced norm $\|\cdot\|$. We denote the sets of complex and extended-complex numbers by $\mathbb{C}$ and $\overline{\mathbb{C}} \coloneq \mathbb{C} \cup \{\infty\}$, respectively. As in \cite{ryu2022scaled}, we avoid $\infty+\infty$, $0/0$, $\infty/\infty$, $0\cdot\infty$ and otherwise adopt the convention $z+\infty=\infty$, $z/\infty=0$, $z/0=\infty$, and $z\cdot\infty=\infty$.  For subsets $S,T\subseteq\overline{\bC}$, set addition is understood as $S+T = \{s+t: s\in S\cap\bC, t\in T\cap\bC\} \cup \{\infty:\ \infty\in S \textnormal{ or } \infty\in T\}$. Scalar multiplication by \(\lambda\in\bC\setminus\{0\}\) is defined by $\lambda S = \{\lambda s: s\in S\cap\bC\}\cup\{\infty: \infty\in S\}$. For $z \in \mathbb{C}$, $z^*$ denotes its complex conjugate whereas, for a bounded, linear operator $M$, $M^*$ denotes its adjoint. We use $\oplus$ to denote the direct sum. For a set-valued mapping $A: \cH \rightrightarrows \cH$, we define its domain $\dom A \coloneq \{x \in \cH \mid A(x) \neq \emptyset\}$, its range $\ran A \coloneq \{u \in \cH \mid u \in A(x) \text{ for some } x \in \cH\}$ and its graph $\graph A \coloneq \{(x, u) \in \mathcal{H}\times\mathcal{H} \mid u \in A(x)\}$. The inverse mapping $A^{-1}:\cH \rightrightarrows \cH$, which always exists, is defined by $A^{-1}(u) \coloneq \{x \in \cH \mid u \in A(x)\}$. The set of zeros is denoted by $\zer(A) \coloneq A^{-1}(0)=\{x \in \cH\mid 0 \in A(x)\}$ and the set of fixed points by $\fix(A) \coloneq \{x \in \cH \mid x \in A(x)\}$. $A$ is firmly nonexpansive if $\langle x - \bar x, u - \bar u \rangle \geq \|u - \bar u\|^2$, $\forall (x,u),(\bar x,\bar u) \in \graph A$, $L$-Lipschitz with $L>0$ if $\|u - \bar u\|\leq L\|x - \bar x\|$, $\forall (x,u),(\bar x,\bar u) \in \graph A$ and contractive if it is $L$-Lipschitz with $L<1$. The composition of two set-valued mappings $A,B : \cH\rightrightarrows \cH$ is defined by $(A \circ B)(x) = \bigcup_{y\in B(x)} A(y)$ for $x\in\cH$. We denote the identity operator by $\id$. The closed disk with center $c\in\mathbb{C}$ and radius $r> 0$ is defined as $D(c, r) \coloneq \{z \in \mathbb{C} : |z - c| \leq r\}$. We denote the half-plane $\mathbb{C}_{\geq \alpha} \coloneq \{z \in \mathbb{C} \mid \realpart(z) \geq \alpha\} \cup \{\infty\}$ with $\alpha \in \mathbb{R}$. Lastly, a set $S \subseteq \overline{\mathbb{C}}$ is said to satisfy the chord property if $z\in S \setminus \{\infty\} \implies \{\theta z +(1-\theta)z^* : \theta\in(0,1)\}\subseteq S$, cf.\ \cite[Fig.\,6]{ryu2022scaled}.

\section{Scaled relative graphs of operator pairs}
We aim to study the relative incremental properties of an operator $A : \cH \rightrightarrows \cH$ compared to another operator $B : \cH \rightrightarrows \cH$. One possible approach is to analyze the composed operator $A \circ B^{-1}$. However, this perspective has several disadvantages: it requires working with explicit operator inverses and compositions, it is more difficult to derive calculus rules, and it may be difficult to preserve problem structure. We therefore take a second approach and study relative incremental properties of operator pairs $(A,B)$.

To this end, we introduce the notion of scaled relative graphs for pairs of operators $(A, B)$, which we show is equivalent to the classical scaled relative graph \cite{ryu2022scaled} of $A \circ B^{-1}$. 
In particular, let $(x, u), (\bar x, \bar u) \in \cH\times\cH$ with $x\neq \bar x$ and define the corresponding complex-conjugate pair
\begin{equation*}
    z_\pm(u - \bar u, x- \bar x) \coloneq \frac{\|u - \bar u\|}{\|x - \bar x\|} \exp(\pm i \angle (x - \bar x, u - \bar u)),
\end{equation*}
where the angle $\angle (x - \bar x, u - \bar u) \in [0, \pi]$ is defined as $\arccos\left(\tfrac{\langle x - \bar x, u - \bar u \rangle}{\|x - \bar x\|\|u - \bar u\|}\right)$ if $x\neq \bar x$ and $u \neq \bar u$, and $0$ otherwise. The SRG of a pair of operators then consists of the union of these pairs, where $A$ and $B$ are evaluated at the same inputs.

\begin{definition}[SRG of a pair of operators] \label{def:pair_srg}
    For $A,B:\cH\rightrightarrows\cH$, the scaled relative graph of a pair of operators $(A, B)$ is
    \[
        \cG(A, B) \coloneq  \left\{z_{\pm}(u_A - \bar u_A, u_B - \bar u_B)\; \Big\vert \; \substack{u_A\in A(x), \bar u_A\in A(\bar x) \\ u_B \in B(x), \bar u_B \in B(\bar x)}\right\}
    \]
    for all $x,\bar x\in\dom A \cap \dom B$ such that $u_B \neq \bar u_B$. Additionally, $\cG(A,B)$ includes $\infty$ if there exist $(x, u_A), (\bar x, \bar u_A) \in \graph A$ and $(x, u_B), (\bar x, \bar u_B) \in \graph B$ such that $u_B = \bar u_B$ (including when $x = \bar x$) and $u_A \neq \bar u_A$. Furthermore, we define the SRG of a class of operator pairs $\cP$ as $\cG(\cP) \coloneq \cup_{(A, B) \in \cP} \cG(A, B)$.
\end{definition}
\begin{remark}
    The classical scaled relative graph \cite{ryu2022scaled} of $A$ is recovered for $B =\id$.
\end{remark}
\begin{remark}
    A related notion has been introduced in \cite[Eq.\,(2.1)]{lee2025convergence}, where $B$ is restricted to a \emph{selection} \cite[p.\,279]{lee2025convergence}.
\end{remark}
We now derive some basic identities that follow from the definition and similar proof techniques as in~\cite{ryu2022scaled}. The main difference arises from $B$ not necessarily being bijective, so handling the $\infty$ case requires extra care.
\begin{proposition}[Basic calculus] \label{thm:scaling_prop}
    Let $A,B,C : \cH \rightrightarrows\cH$, and $\alpha,\beta \in \bR\setminus\{0\}$. Then,
    \begin{propenum}
        \item $\cG(A,B) = \cG(A \circ B^{-1}, \id)$. \label{thm:equiv}
        \item $\cG(\alpha A, \beta B) = (\alpha/\beta)\cG(A, B).$ \label{thm:scaling_prop:scale}
        \item $\cG(B, A) = (\cG(A, B))^{-1} \coloneq \{(z^{-1})^*\mid z \in \cG(A, B)\}$.\label{thm:scaling_prop:inv}
        \item $\cG( \id, A) = (\cG(A, \id))^{-1} = \cG(A^{-1}, \id)$. \label{thm:scaling_prop:id}
        \item If either $\cG(A,C)$ or $\cG(B,C)$ satisfies the chord property, then $\cG(A + B, C) \subseteq \cG(A, C) + \cG(B, C)$. \label{thm:scaling_prop:trans}
        \item If either $\cG(B, A)$ or $\cG(C, A)$ satisfies the chord property, then $\cG(A, B + C) \subseteq (\cG(A, B)^{-1} + \cG(A, C)^{-1})^{-1}$.\label{thm:scaling_prop:trans_inv}
    \end{propenum}
\end{proposition}

If one of the operators is single-valued or satisfies some mild additional properties, even more calculus rules can be established.
\begin{proposition}[Additional calculus] \label{thm:single_valued}
    Let $A,B : \cH\rightrightarrows\cH$ and suppose $\sing: \cH \to \cH$ is single-valued. Then,
    \begin{propenum}
        \item If $\sing$ is not constant, then $\cG(\sing,\sing) = \{1\}$.\label{thm:single_valued:id}
        \item $\cG(A + \sing, \sing) = 1 + \cG(A, \sing)$.\label{thm:single_valued:sum}
        \item $\cG (\sing \circ (A + \sing)^{-1}, \id) = \cG(\sing, A + \sing)$.  \label{thm:single_valued:sum_inv}
        \item $\cG(A\circ \sing, B \circ \sing) \subseteq \cG(A, B)$ with equality if $\sing$ is surjective. \label{thm:single_valued:precomp}
        \item If $\cG(\sing, \id) \subseteq D(0, L)$ and $\cG(A,B) \subseteq D(0,l)$ for some $L,l>0$, then $\cG(\sing \circ A, B) \subseteq D(0, Ll)$. \label{thm:single_valued:lip}
        \item If $\cG(A, \sing) \subseteq D(0, L)$ and $\cG(\sing, \id) \subseteq D(0, l)$ for some $L,l>0$, then $\cG(A, \id) \subseteq D(0,Ll)$. \label{thm:single_valued:lip2}
        \item If $M$ is a bounded, invertible linear operator on $\cH$ and $\cG(A, B) \subseteq \bC_{\geq 0}$, then $\cG((M^{-1})^* \circ A, M \circ B) \subseteq \bC_{\geq 0}$. \label{thm:single_valued:linmap}
    \end{propenum}
\end{proposition}
\begin{remark}
    Let $A : \mathcal{H}\rightrightarrows \mathcal{H}$ and $F :\mathcal{H}\to \mathcal{H}$. We have by \Cref{thm:equiv} that $\mathcal{G}(A, F) = \mathcal{G}(A \circ F^{-1}, \id)$. In particular, consider $\mathcal{H} = \mathbb{R}^n$ and a positive definite matrix $W \in \mathbb{R}^{n\times n}$, which admits a decomposition $W = X^\top X$ with $X \in \mathbb{R}^{n\times n}$ nonsingular by \cite[Thm.\,7.2.7]{horn2012matrix}. Then, $\mathcal{G}(X \circ A, X) = \mathcal{G}(X\circ A \circ X^{-1}, \id)$ corresponds to the incremental variant of the weighted scaled graph from \cite{de2025exploiting}.
\end{remark}

A crucial property in the context of the original SRG is the concept of SRG-fullness of operator classes \cite[Sec.\,3.3]{ryu2022scaled}, since these allow for membership checking based on geometric containment of the SRG, which forms the final step of the approach shown in \Cref{fig:meta_proc}. A generalization to our framework is straightforward.
\begin{definition}[SRG-full operator classes of pairs] \label{def:srgfull}
A class $\cP$ of operator pairs is SRG-full if \[(A, B) \in \mathcal{P} \iff \cG(A, B) \subseteq \cG(\cP).\]
\end{definition}

Similar to \cite[Thm.\,2]{ryu2022scaled}, classes defined through some nonnegatively homogeneous function $h : \bR^3\to\bR$, i.e., $h(\kappa a, \kappa b, \kappa c) = \kappa h(a, b, c)$ for all $\kappa\geq0$, satisfy this desirable property. In fact, if $\cP$ is SRG-full, then such a nonnegatively homogeneous function always exists.
\begin{proposition} \label{prop:srgfull}
    Let $\cP$ be a class of operator pairs. Then $\cP$ is SRG-full if and only if there exists some nonnegatively homogeneous function $h:\bR^3\to \bR$ such that $(A, B) \in \cP$ if and only if $\forall (x, u_A), (\bar x, \bar u_A) \in \graph A, (x, u_B), (\bar x, \bar u_B) \in \graph B$:
    \begin{equation} \label{eq:srgfull_hom}
         h(\|u_A - \bar u_A\|^2, \|u_B - \bar u_B\|^2, \langle u_A - \bar u_A, u_B - \bar u_B\rangle) \leq  0.
    \end{equation}
\end{proposition}
\begin{remark}
    Note that the choice of $h$ representing a class is not unique, and that there may exist representations that are not nonnegatively homogeneous, e.g., $h_1(a,b,c)=-c$, $h_2(a,b,c) = -2c$, and $h_3(a,b,c)=-c^3$ all represent (pair of) monotone operators, but $h_3$ clearly violates the preceding homogeneity condition.
\end{remark}
In the context of the original SRG, one particularly interesting SRG-full class is the one defined by $h : (a, b, c) \mapsto \rho a + \mu b - c$, which covers $(\mu, \rho)$-semimonotone operators, see \cite[Prop.\,3.2]{quan2024scaled}, a class that was introduced in \cite[Def.\,4.1]{evens2025convergence}. In the following, we generalize this class to the setting of operator pairs, inspired by the operator-pair monotonicity introduced in \cite[Def.\,3]{borwein2011characterization}, \cite[Eq.\,(5)]{le2025solving}. Note that classical $(\mu,\rho)$-semimonotonicity is recovered by taking $B=\id$.
\begin{definition}[Semimonotone operator pairs] \label{def:semimon}
     Let $A,B:\cH \rightrightarrows \cH$ and $\mu,\rho\in\bR$. Then $(A, B)$ is $(\mu, \rho)$-semimonotone if $\forall (x, u_A), (\bar x, \bar u_A) \in \graph A, (x, u_B), (\bar x, \bar u_B) \in \graph B:$
    $
        \langle u_A - \bar u_A, u_B - \bar u_B\rangle \geq \mu\|u_B - \bar u_B\|^2 + \rho \|u_A - \bar u_A\|^2.
    $
     The class of all $(\mu,\rho)$-semimonotone operator pairs is denoted by $\cS_{\mu,\rho}$.
\end{definition}
We also denote $\cM_\mu \coloneq \cS_{\mu,0}$ for the $\mu$-monotone pairs of operators. When $B=\id$ and $\mu > 0$, we recover the class of $\mu$-strongly monotone operators, while we recover the class of $|\mu|$-hypomonotone operators when $\mu < 0$, similarly to \cite[Rem.\,4.2]{evens2025convergence}. The SRG of the class of $(\mu,\rho)$-semimonotone operator pairs is now shown to be exactly the same as in \cite[Prop.\,3.4]{quan2024scaled}, for which we moreover derive an alternative representation.
\begin{proposition}[SRG of semimonotone operator pairs] \label{prop:srg_semi}
    Let $\mu, \rho \in \bR$. Then, $\cS_{\mu,\rho}$ is SRG-full. Moreover,
    \[
        \cG(\cS_{\mu,\rho}) = \{z \in \bC \mid \realpart(z) \geq \mu + \rho|z|^2\} \, (\cup \{\infty\} \textnormal{ if $\rho\leq 0$)}.
    \]
    In particular,
    \begin{equation} \label{eq:srg_mon}
        \cG(\cM_\mu) = \{z\in\bC \mid \realpart(z) \geq \mu\} \cup \{\infty\}.
    \end{equation}
\end{proposition}
\begin{iproof}
    The forward inclusion is shown using the definition of semimonotone operator pairs, while the reverse inclusion follows because classical semimonotone operators (i.e., with $B=\id$) are contained in $\cS_{\mu,\rho}$, so that \cite[Prop.\,3.4]{quan2024scaled} applies.
\end{iproof}

\begin{example} \label{ex:srg-linear}
    Consider the linear operator $A_{\rm lin} = \begin{bsmallmatrix}
        1/2 & 2 \\ -1/2 & 1/2
    \end{bsmallmatrix} \oplus [2]$ from \cite[Fig.\,5]{ryu2022scaled}. \Cref{fig:srg-linear} visualizes, for different operators \(B\), the numerical SRG of the pair \((A_{\rm lin}, B)\).
    Here and henceforth, numerical SRGs are computed with each coordinate sampled independently and uniformly from $[-1000, 1000]$.
    For \(B = \id\), this reduces to the standard SRG of \(A_{\rm lin}\). 
    Inspired by \cite[Lem.\,5.1]{le2025solving}, we also consider the choice $B = A_{\rm lin} + 2 \kappa \id$, where \(\kappa = 0.45 \vert \alpha \vert\) and \(\vert \alpha \vert\) is the smallest eigenvalue of \(A_{\rm lin}\) in absolute value.
    Finally, we also consider \(B = 2A_{\rm lin}^{-\top}\).

    It is clear that $(A_{\rm lin}, \id)$ is not monotone by \eqref{eq:srg_mon}. Further, the second and third SRGs suggest that $(A_{\rm lin}, A_{\rm lin} + 2\kappa\id), (A_{\rm lin}, 2A_{\rm lin}^{-\top}) \in \cM_0$, though this cannot be concluded by sampling alone. That they indeed have this property follows from \cite[Lem.\,5.1]{le2025solving} and \Cref{thm:linmon:fullrank}, respectively.
    \begin{figure}[hbpt]
        \centering
        \subfloat[$B = \id$]{
            \includegraphics[page=1,width=0.14\textwidth]{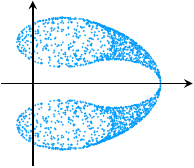}
        } \hfill
        \subfloat[$B = A_{\rm lin} + 2 \kappa \id$]{
            \includegraphics[page=2,width=0.14\textwidth]{figures/srg_ryu.pdf}
        } \hfill
        \subfloat[$B = 2 A_{\rm lin}^{-\top}$]{
            \includegraphics[page=3,width=0.14\textwidth]{figures/srg_ryu.pdf}
        }
        \caption{Numerical SRG of \((A_{\rm lin}, B)\).} \label{fig:srg-linear}
    \end{figure}
\end{example}
\begin{example} \label{ex:srg-transistor}
      Let $A_{\rm NPN}$ be the NPN transistor modeled by the Ebers--Moll model, see \cite[Sec.\,4.2]{quan2024scaled},
    \[
        A_{\rm NPN} \begin{bmatrix}
            v_1 \\ v_2
        \end{bmatrix} \coloneq \left\{R\begin{bmatrix}
            u_1 \\ u_2
        \end{bmatrix} \mathrel{\bigg|}
        \begin{aligned}
            &u_1 {}\in{} A_{\rm D}(v_1)\\
            &u_2 {}\in{} A_{\rm D}(v_2) 
        \end{aligned} \right\}
    \]
    where $R\coloneq\left[\begin{smallmatrix}
        1 & -\alpha_R \\ -\alpha_F & 1
    \end{smallmatrix}\right]$, $0 \leq \alpha_R, \alpha_F< 1$, and $A_{\rm D}$ is a diode model satisfying $(A_{\rm D},\id)\in \cM_0$. 
    \Cref{fig:srg-transistor} shows, for different operators $B$, the numerical SRGs of the pairs \((A_{\rm NPN}, B)\).
    For \(B = \id\), we observe that the NPN transistor is angle-bounded \cite[Def.\,3.6]{quan2024scaled}, as proven in \cite[Prop.\,4.4]{quan2024scaled}, and that it is not monotone by \eqref{eq:srg_mon}.
    Unlike in \cref{ex:srg-linear}, the choice \(B = A_{\rm NPN} + 2 \kappa \id\), where \(\kappa = 0.45 \vert \alpha \vert\) and \(\vert \alpha \vert\) is the smallest eigenvalue of \(R\) in absolute value \cite[Lem.\,5.1]{le2025solving}, does not yield a monotone pair of operators.
    Finally, for \(B = \det (R) R^{-\top}\), the sampled SRG suggests monotonicity of the pair \((A_{\rm NPN}, \det(R) R^{-\top})\), which \Cref{cor:trans} later establishes formally.
    \begin{figure}[t]
        \centering
        \subfloat[$B = \id$]{%
            \includegraphics[page=1,width=0.15\textwidth]{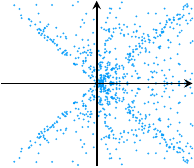}
        } \hfill
        \subfloat[$B = A_{\rm NPN} + 2 \kappa \id$]{
            \includegraphics[page=2,width=0.15\textwidth]{figures/srg_transistor.pdf}
        } \hfill
        \subfloat[$B = \det(R) R^{-\top}$]{
            \includegraphics[page=3,width=0.15\textwidth]{figures/srg_transistor.pdf}
        }
        \caption{Numerical SRG of \((A_{\rm NPN}, B)\).} \label{fig:srg-transistor}
    \end{figure}  
\end{example}

\section{Properties of nonlinear resolvents}
Let $A:\cH\rightrightarrows\cH$ and consider the zero inclusion problem 
\[
    0 \in A(x)
\]
that arises ubiquitously in optimization and systems theory. A classical way to solve this problem is to pose it as finding a fixed point of a related operator. The most well-known example is the \emph{resolvent} $J_{\gamma A} = (\gamma A + \id)^{-1}$, which leads to the celebrated proximal point algorithm. Often, this fixed point operator is shown to be firmly nonexpansive or contractive, from which convergence readily follows by the \KM\ or Banach fixed-point theorem.

We now recall two nonlinear resolvents called the warped resolvent \cite{bui2020warped} and the transformed resolvent~\cite{le2025solving}.
\begin{definition}[Nonlinear resolvents]
    Let $A: \cH\rightrightarrows\cH$ and suppose $\sing:\cH\to\cH$ is single-valued. The transformed resolvent of $A$ with respect to $\sing$ is defined as $T_A^\sing \coloneq \sing \circ (A + \sing)^{-1}$. The warped resolvent of $A$ with respect to $\sing$ is defined as $J_A^\sing \coloneq (A + \sing)^{-1}\circ \sing$, provided that $\ran \sing \subseteq \ran(A + \sing)$.
\end{definition}
These resolvents are useful since $\fix T_{\gamma A}^\sing = \sing(\zer A)$ and $\fix J_{\gamma A}^\sing = \zer A$ in light of \cite[Prop.\,3.2]{le2025solving}. Therefore, if we can show firm nonexpansiveness or contraction, then standard theory shows convergence, as described above. These properties have indeed been shown under the pair of monotonicity framework in \cite{le2025solving}. For the remainder of this paper, we assume that the considered resolvents are everywhere defined on the relevant space.

We now provide the geometric picture of these derivations by giving purely SRG-based proofs. This approach yields concise proofs that capture the core insights and also allow us to see that the established results are tight.

\begin{figure}[t]
    \centering
    \subfloat{%
        \includegraphics[page=1,width=0.15\textwidth]{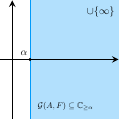}
    } \hfill
    \subfloat{
        \includegraphics[page=2,width=0.15\textwidth]{figures/srgv2.pdf}
    } \hfill
    \subfloat{
        \includegraphics[page=3,width=0.15\textwidth]{figures/srgv2.pdf}
    }
    \caption{Geometry of proof \Cref{prop:transf_resv}.} \label{fig:proof_transf}
\end{figure}

First, we show that if a pair of operators is $\alpha$-monotone for some $\alpha\geq 0$, then the transformed resolvent is firmly nonexpansive or even contractive. We then show a similar result for the warped resolvent, under some additional assumptions on $F$.
\begin{proposition}[Properties transformed resolvent] \label{prop:transf_resv}
    Let $A: \cH\rightrightarrows\cH$, $\gamma>0$, and suppose $\sing:\cH\to\cH$ is single-valued. If $(A, \sing)\in\cM_\alpha$ with $\alpha\geq0$, then the transformed resolvent $T_{\gamma A}^\sing$ has $\cG(T_{\gamma A}^\sing, \id) \subseteq D(1/(2+2\gamma\alpha), 1/(2+2\gamma\alpha))$. In particular, if $\alpha=0$, then the transformed resolvent is firmly nonexpansive. If $\alpha>0$, then the transformed resolvent is Lipschitz continuous with constant $1/(1+\alpha\gamma) < 1$, i.e., contractive.
\end{proposition}
\begin{proof} The geometry of this proof is shown in \Cref{fig:proof_transf}.

    Since $(A, \sing)\in\cM_\alpha$, we find from \eqref{eq:srg_mon} that $\cG(A, \sing) \subseteq \bC_{\geq\alpha}$. Then, by \Cref{thm:scaling_prop:scale}, we obtain $\cG(\gamma A, \sing) \subseteq \bC_{\geq\gamma\alpha}$. Further, from the single-valuedness of $\sing$, \Cref{thm:single_valued:sum} ensures that $\cG(\gamma A + \sing, \sing) \subseteq 1 + \bC_{\geq \gamma\alpha}$. Lastly, the property \Cref{thm:single_valued:sum_inv} and the inversion rule \Cref{thm:scaling_prop:inv} yield that $\cG(T_{\gamma A}^\sing, \id) = \cG(\sing, \gamma A + \sing) \subseteq (1 + \bC_{\geq \gamma \alpha})^{-1} = D(1/(2+2\gamma\alpha), 1/(2+2\gamma\alpha))$.

    Since $\cG(T_{\gamma A}^\sing, \id)$ is the standard SRG of $T_{\gamma A}^\sing$, it follows from \cite[Prop.\,1 and Thm.\,2]{ryu2022scaled} that $T_{\gamma A}^\sing$ is firmly nonexpansive if $\alpha=0$ and Lipschitz continuous with factor $1/(1+\alpha\gamma) < 1$ if $\alpha>0$.
\end{proof}

\begin{proposition}[Lipschitz continuity warped resolvent] \label{prop:warp_resv}
    Let $A: \cH\rightrightarrows\cH$, $\gamma>0$, and suppose $\sing: \cH\to\cH$ is single-valued. If $(A, \sing)\in \cM_\alpha$ with $\alpha\geq0$, then the warped resolvent satisfies $\cG(\sing\circ J_{\gamma A}^\sing, \sing) \subseteq D(1/(2+2\gamma \alpha), 1/(2+2\gamma \alpha))$. If, moreover $\sing^{-1}$ is well-defined and $l$-Lipschitz, and $\sing$ is $L$-Lipschitz, then $\cG(J_{\gamma A}^\sing, \id) \subseteq D(0,Ll/(1+\gamma\alpha))$.
\end{proposition}

\begin{figure}[t]
    \centering
    \subfloat{%
        \includegraphics[page=4,width=0.15\textwidth]{figures/srgv2.pdf}
    } \hfill
    \subfloat{
        \includegraphics[page=5,width=0.15\textwidth]{figures/srgv2.pdf}
    } \hfill
    \subfloat{
        \includegraphics[width=0.15\textwidth]{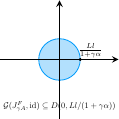}
    }
    \caption{Geometry of proof \Cref{prop:warp_resv}.} \label{fig:proof_warp}
\end{figure}
\begin{proof}
    The geometry of the proof is shown in \Cref{fig:proof_warp}.
    
    Since $\sing\circ J_{\gamma A}^\sing = T_{\gamma A}^\sing \circ \sing$, \Cref{prop:transf_resv} and the precomposition rule \Cref{thm:single_valued:precomp} ensure that $\cG(\sing\circ J_{\gamma A}^\sing, \sing) \subseteq D(1/(2+2\gamma \alpha), 1/(2+2\gamma\alpha))$. If, moreover $\sing^{-1}$ is well-defined and $l$-Lipschitz, then $\cG(\sing^{-1}, \id) \subseteq D(0,l)$ by \cite[Prop.\,1]{ryu2022scaled}. From \Cref{thm:single_valued:lip} and the fact that $D(1/(2+2\gamma \alpha), 1/(2+2\gamma\alpha))\subseteq D(0,1/(1+\gamma\alpha))$, it then follows that $\cG(\sing^{-1} \circ \sing \circ J_{\gamma A}^\sing, \sing) = \cG(J_{\gamma A}^\sing, \sing) \subseteq lD(0, 1/(1+\gamma\alpha)) = D(0, l/(1+\gamma\alpha))$. Lastly, by \Cref{thm:single_valued:lip2} and the fact that $\cG(\sing, \id) \subseteq D(0,L)$ from \cite[Prop.\,1]{ryu2022scaled}, we conclude that $\cG(J_{\gamma A}^\sing, \id) \subseteq D(0,Ll/(1+\gamma\alpha))$.
\end{proof}

\begin{remark} \label{rem:str_mon_diff}
    In \cite{le2025solving}, strong monotonicity of pairs of operators is defined differently from our $\cM_\alpha$. Nevertheless, if $(A, B)$ satisfies \cite[Eq.\,(6)]{le2025solving}, i.e., $\forall (x,u_A),(\bar x, \bar u_A) \in \graph A$, $(x,u_B),(\bar x,\bar u_B)\in\graph B$, \[\langle u_A - \bar u_A, u_B - \bar u_B\rangle \geq \alpha \|x - \bar x\|^2\] and $B$ is $L$-Lipschitz continuous, then $\|x - \bar x\|^2 \geq L^{-2}\|u_B - \bar u_B\|^2$ and it follows that $(A, B) \in \cM_{\alpha L^{-2}}$ by \Cref{def:semimon}, so \Cref{prop:transf_resv} exactly recovers \cite[Prop.\,3.4]{le2025solving}.
\end{remark}
\begin{example} \label{ex:forward}
    Note that the definition of semimonotone operator pairs allows for great flexibility in the choice of $(A, F)$. For $\mu = \rho = 0$ and $F = \id - A$, \cref{def:semimon} recovers the class of firmly nonexpansive operators:
    \begin{align*}
        &\langle A(x)-A(\bar x),x-A(x)-(\bar x-A(\bar x)) \rangle \geq 0
        \\
        \iff &\langle A(x)-A(\bar x),x-\bar x \rangle \geq \|A(x)-A(\bar x)\|^2,
    \end{align*}
    for all $x, \bar x \in \dom A$. In that case, the transformed resolvent becomes $T_A^F = F = \id - A$, i.e., the standard forward step. 
\end{example}
        
\begin{example} 
    We can recover more relaxed conditions by using the so-called nonlinear preconditioning technique \cite{maddison2021dual,laude2025anisotropic,laude2023anisotropic,oikonomidis2025nonlinearlypreconditionedgradientmethods}: choosing $F = \id - \nabla \psi \circ A$ where $\nabla \psi$ is the gradient of a Legendre function \cite[p.\,6]{laude2025anisotropic} and $\ran A \subseteq \dom \nabla \psi$. In that case, the semimonotonicity inequality with $\mu=\rho=0$ takes the form $
        \langle A(x)-A(\bar x),x-\bar x \rangle 
   \geq
        \langle \nabla \psi (A(x))-\nabla \psi(A(\bar x)), A(x)-A(\bar x) \rangle.
    $
    
    Clearly, this implies that $A$ is a monotone operator while if $\psi = \tfrac{1}{2}\|\cdot\|^2$, we recover \Cref{ex:forward}. We remark that by choosing a suitable $\psi$, we can make the inequality above less restrictive than the one in \Cref{ex:forward} as shown in \cite{oikonomidis2025nonlinearlypreconditionedgradientmethods}. The corresponding transformed resolvent becomes $T_A^\sing = (\id - \nabla \psi \circ A)\circ ((\id - \nabla \psi) \circ A + \id)^{-1}$.
\end{example}
\begin{example} 
    As in \cref{ex:forward}, let \(F = \id - A\). 
    Suppose that \(A = \gamma \nabla \psi \circ \nabla f\), where \(f : \bR^2 \to \bR : x \mapsto \frac{1}{4}\sum_{i = 1}^2 x_i^4\) and \(\gamma = 0.1\) is a fixed step size.
    Then, the transformed resolvent reduces to a nonlinearly preconditioned forward step \(T_A^F = \id - \gamma \nabla \psi \circ \nabla f\) \cite{oikonomidis2025nonlinearlypreconditionedgradientmethods}. Typically, $\nabla \psi$ satisfies $\nabla\psi(y)=0$ if and only if $y=0$, in which case the zeros of $A$ correspond to the stationary points of $f$.
    
    \Cref{fig:srg-aniso} visualizes the numerical SRG of the pair \((A, F)\) for different separable nonlinear preconditioners of the form \(\nabla \psi = (h'(y_1), h'(y_2))\). 
    Without preconditioning, i.e., for \(h'(y) = y\), the pair is clearly not monotone. 
    For hard clipping \(h'(y) = \proj_{[-1, 1]}(y)\), and for \(h'(y) = \arcsinh(y)\) the numerical SRG suggests monotonicity of the pair \((A, F)\).
    \begin{figure}[t]
        \centering
        \subfloat[$h'(y) = y$]{%
            \includegraphics[page=1,width=0.15\textwidth]{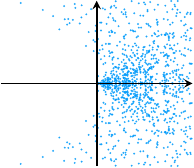}
        } \hfill
        \subfloat[$h'(y) = \proj(y)$]{
            \includegraphics[page=2,width=0.15\textwidth]{figures/srg_aniso.pdf}
        } \hfill
        \subfloat[$h'(y) = \arcsinh(y)$]{
            \includegraphics[page=3,width=0.15\textwidth]{figures/srg_aniso.pdf}
        }
        \caption{Numerical SRG of \((\gamma \nabla \psi \circ \nabla f, \id - \gamma \nabla \psi \circ \nabla f)\) for \(f(x) = \frac{1}{4} \sum_{i=1}^2 x_i^4\) and \(\nabla \psi = (h'(y_1), h'(y_2))\).
        We denote by \(\proj\) the Euclidean projection onto the interval \([-1, 1]\).
        } \label{fig:srg-aniso}
    \end{figure}
\end{example}

Lastly, inspired by \cite[Ex.\,2.3 and 2.4]{le2025solving}, we provide two more examples to showcase the SRG approach.

\begin{example}
    Let $A : \cH \rightrightarrows \cH$. Suppose that $A = B + \sing$, where $B:\cH\rightrightarrows\cH$ and $\sing : \cH\to\cH$ is single-valued. If $(\sing, B) \in \cM_0$, then $(A, \sing) \in \cM_1$.
\end{example}
\begin{proof}
    Since $(\sing, B) \in \cM_0$, we have that $\cG(\sing, B) \subseteq \bC_{\geq 0}$. Then, from the inversion rule \Cref{thm:scaling_prop:inv}, $\cG(B, \sing) \subseteq (\bC_{\geq 0})^{-1} = \bC_{\geq 0}$. Further, by \Cref{thm:single_valued:sum}, $\cG(B + \sing, \sing) \subseteq 1+ \bC_{\geq 0} = \bC_{\geq 1}$ and we conclude that $(A, \sing) = (B + \sing, \sing) \in \cM_1$ by \Cref{prop:srg_semi}.
\end{proof}

\begin{example}  \label{ex:sum}
    Let $A,B:\cH\rightrightarrows\cH$ and suppose $\sing:\cH\to\cH$ is single-valued. If $(A, \id)\in \cM_0$ and $(B, \sing)\in \cM_\alpha$ for some $\alpha\geq0$, then $(A \circ \sing + B, \sing)
    \in \cM_\alpha$.
\end{example}
\begin{proof}
     Since $(A, \id) \in \cM_0$ and $(B,\sing) \in \cM_\alpha$, we have $\cG(A, \id) \subseteq \bC_{\geq 0}$ and $\cG(B, \sing) \subseteq \bC_{\geq\alpha}$. Then, by the precomposition rule, \Cref{thm:single_valued:precomp}, $\cG(A \circ \sing, \sing) \subseteq \bC_{\geq 0}$, and from \Cref{thm:scaling_prop:trans} (passing to chord completions if necessary, see \cite[Def.\,4]{krebbekx2025graphical}), we obtain that $\cG(A\circ \sing + B, \sing) \subseteq \bC_{\geq0} + \bC_{\geq\alpha} = \bC_{\geq\alpha}$, which readily leads to the desired result by \Cref{prop:srg_semi}.
\end{proof}

Lastly, we provide a result that recovers part of \cite[Cor.\,25.6]{bauschke2017correction} when $\sing=\id$, and is useful for deriving (linearly) preconditioned algorithms, as we will do in \Cref{ex:circuit}.
\begin{proposition} \label{prop:genprox}
    Let $A:\cH\rightrightarrows\cH$ and suppose $\sing : \cH \to \cH$ is single-valued. Let $M$ be a bounded, invertible linear operator on $\cH$. If $(A, \sing) \in \cM_0$, then $((M^{-1})^*AM^{-1}, M\sing M^{-1}) \in \cM_0$.
\end{proposition}
\begin{proof}
    Since $(A, \sing)\in \cM_0$, we have $\cG(A, \sing) \subseteq \bC_{\geq0}$. Then, from \Cref{thm:single_valued:linmap}, we obtain $\cG((M^{-1})^*A, M\sing) \subseteq \bC_{\geq0}$ and by precomposition, \Cref{thm:single_valued:precomp}, we have $\cG((M^{-1})^*A M^{-1}, M\sing M^{-1}) \subseteq \bC_{\geq0}$, so we conclude that $((M^{-1})^* A M^{-1}, M\sing M^{-1}) \in \cM_0$ by \Cref{prop:srg_semi}.
\end{proof}

\section{Application to circuit theory}
In this section, we apply the paired monotonicity framework to solve two inclusions involving nonsmooth, multivalued and highly nonmonotone operators in the context of circuit theory. To this end, we first show that linear mappings composed with monotone mappings fit naturally into this framework, a class that includes the NPN transistor. The SRG framework can be used as a simple visual tool for detecting when paired monotonicity fails and may also certify tightness when monotonicity does hold (see \Cref{fig:srg-transistor}).

\begin{figure}[t]
    \centering
    \subfloat[Leakage transistor]{%
        \begin{circuitikz}[american]
            \draw (0,0) node[npn, rotate=90] (p) {};
            \draw (p.B) to[short, -*] (0, -2) coordinate (B) to (B -| p.C) to[short, i^<=$i_1$] ++(-0.6, 0) coordinate (A) to[open, o-o, v^>=$v_1$] (A |- p.C) to[short, i<=$i_1$] (p.C);
            \draw (B) to (B -| p.E) to[short, i_<=$i_2$] ++(0.6,0) coordinate (C) to[open, o-o, v^>=$v_2$] (C |- p.C) to[short,i_<=$i_2$] (p.E);
            \draw (p.E) to[/tikz/circuitikz/bipoles/length=1.25cm, R, l_=$r$] (p.E |- C);
            \draw (p.C) to[/tikz/circuitikz/bipoles/length=1.25cm, R, l^=$r$] (p.C |- A);
        \end{circuitikz}
        \label{fig:NPN-leaky}
    } \hfill
    \subfloat[Common-emitter amplifier]{%
        \begin{circuitikz}[american]
            \ctikzset{bipoles/resistor/height=0.1}
            \ctikzset{bipoles/resistor/width=0.3}
            \ctikzset{bipoles/ageneric/height=0.2}
            \ctikzset{bipoles/ageneric/width=0.45}
    
            \draw (0,0) node[npn] (p) {};
            \draw (p.E) to[R, -, l^=$R_E$] ++(0, -0.65) to[short, -] ++(0, -0.2) coordinate (G) to[short] ++(-1.35, 0) coordinate (A) to[/tikz/circuitikz/bipoles/length=1cm, V, invert, l=$v_{\textnormal{in}}$] (A |- p.B) to (p.B);
        
            \draw (p.C) to[R, l_=$R_C$] ++(1.5, 0) coordinate (B) to[/tikz/circuitikz/bipoles/length=1cm, V, l= $v_+$] (B |- A) to ++(A |- p.B) node[ground]{} to (G);

            \node at (0.3, -0.7) {$+$};
            \node at (0.3, -1.45) {$-$};

            \node at (0.3, 0.45) {$+$};
            \node at (1.2, 0.45) {$-$};            
        \end{circuitikz}
        \label{fig:common-emitter-amplifier}
    }
    \caption{Examples of nonmonotone electrical circuits.}
\end{figure}
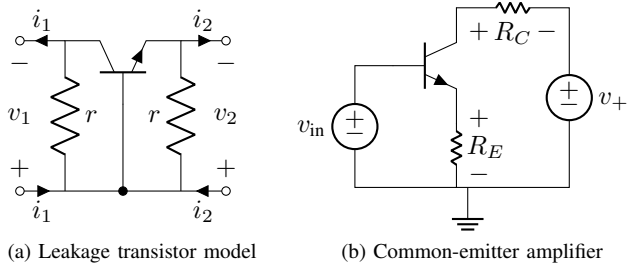

\begin{proposition}
    Let $M \in \bR^{n\times n}$ and $B:\bR^n \rightrightarrows \bR^n$. Let $A = M \circ B$. Suppose that $(B, \id)\in \cM_0$.
    \begin{propenum}
        \item If $M$ is nonsingular, then $(A, cM^{-\top}) \in \cM_0$ for any $c>0$. \label{thm:linmon:fullrank}
        \item If $\rank(M) = n-1$, then $(A, yx^\top) \in \cM_0$ for some $x,y\in\bR^n\setminus\{0\}$. \label{thm:linmon:rankdef}
    \end{propenum}
\end{proposition}
\begin{proof}
    ``\labelcref{thm:linmon:fullrank}'': Let $(x, u), (\bar x, \bar u) \in \graph(cM^{-\top})$ and let $(x, u_A), (\bar x, \bar u_A) \in \graph A$. Then,
    \begin{align*}
        \langle u - \bar u, u_A - \bar u_A\rangle &= \langle cM^{-\top}(x - \bar x), M(u_B - \bar u_B)\rangle\\
        &= c\langle x - \bar x, u_B - \bar u_B\rangle
    \end{align*}
    for $u_B \in B(x), \bar u_B\in B(\bar x)$. From the monotonicity of $B$, we find that $\langle x - \bar x, u_B - \bar u_B\rangle \geq 0$, so $(A, cM^{-\top}) \in \cM_0$.

    ``\labelcref{thm:linmon:rankdef}'': By a similar proof as in \Cref{thm:linmon:fullrank}, $(A, \adj M^\top) \in \cM_0$ provided that $(\det M) \geq 0$, where $\adj M$ denotes the adjugate of $M$, i.e., the transpose of the cofactor matrix. To derive this, the defining property of the adjugate $(\adj M)M = (\det M)I$ is used. Further, if $\rank(M)=n-1$, then $\det(M)=0$ and $\rank(\adj M)=1$, so $\adj M$ admits the factorization $\adj M = xy^\top$ for some $x,y\in\bR^n\setminus\{0\}$ and the result is proven. See \cite[Sec.\,0.8.2]{horn2012matrix} for a more detailed discussion on the used properties of adjugates.
\end{proof}

\begin{corollary} \label{cor:trans}
    The operator $A_{\rm NPN}$ in \Cref{ex:srg-transistor} satisfies $(A_{\rm NPN}, B) \in \cM_0$, where $B=(\det R)R^{-\top} =\begin{bsmallmatrix}
        1&\alpha_F \\ \alpha_R & 1
    \end{bsmallmatrix}$.
\end{corollary}

\begin{figure}[t]
    \centering \vspace{0.15cm}
    \includegraphics[width=0.48\linewidth]{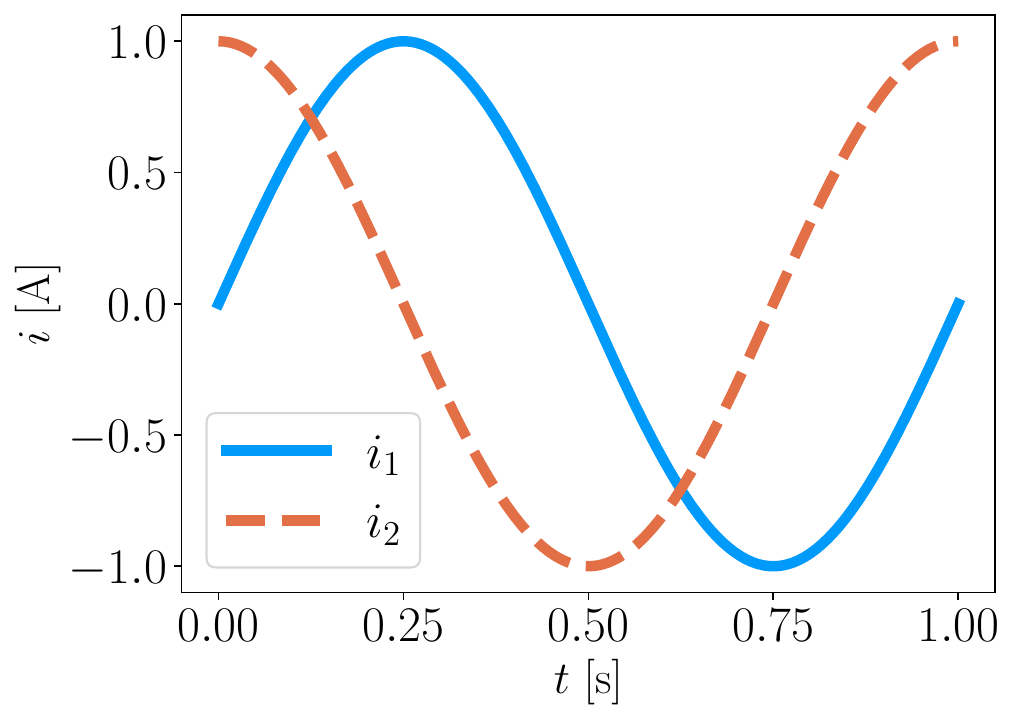}
    \hfill
    \includegraphics[width=0.48\linewidth]{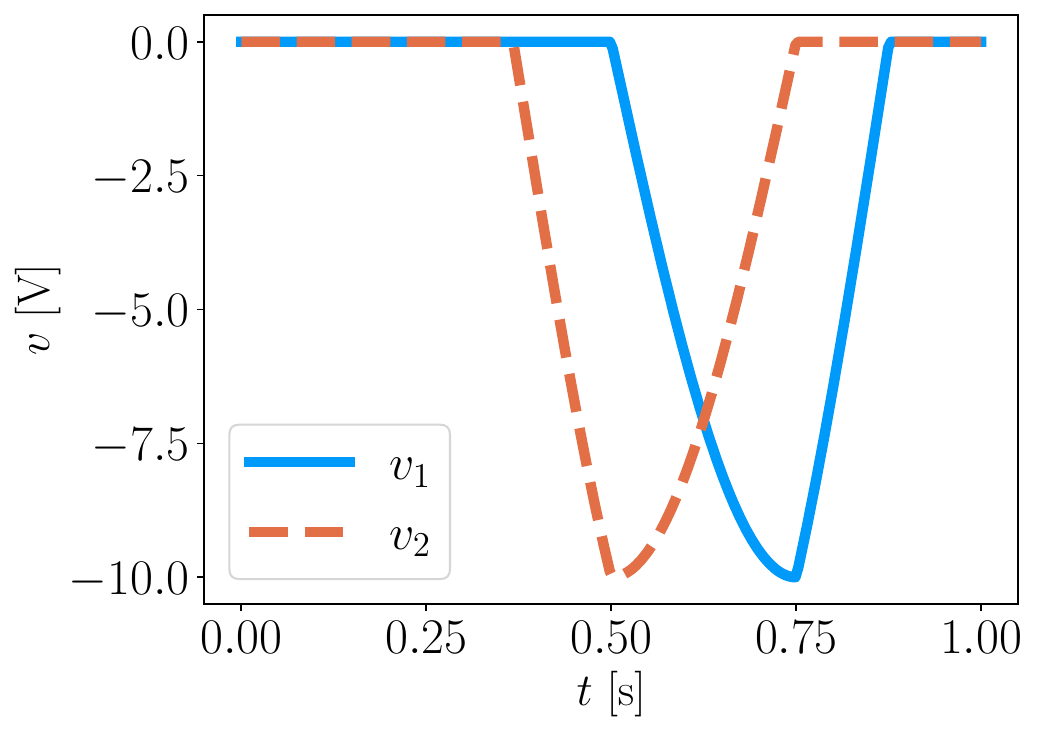}
    \caption{Solution to the inclusion problem \eqref{eq:transistor_eq} for the same configuration as in \cite[Fig.\,6]{quan2024scaled}: leakage resistance \(r = \SI{10}{\ohm}\) and a desired sinusoidal current \(i\) are given.} \label{fig:leaky_trans}
\end{figure}

In the following example\footnote{The code for reproducing the experiments can be found at \url{https://github.com/alexanderbodard/SRGs_for_pairs}.}, we consider the same experiment as in \cite[Prop.\,5.1]{quan2024scaled} and provide a transformed proximal point iteration that converges without any step size restriction.
\begin{example}
    We first consider the leaky transistor shown in \Cref{fig:NPN-leaky}. The associated inclusion problem is
    \begin{equation} \label{eq:transistor_eq}
        \begin{bmatrix}
            i_1 \\ i_2
        \end{bmatrix} \in A_{\rm NPN} \begin{bmatrix}
            v_1 \\ v_2
        \end{bmatrix} + \frac{1}{r}\begin{bmatrix}
            v_1 \\ v_2
        \end{bmatrix} \eqcolon A_{{\rm NPN}, r} \begin{bmatrix}
            v_1 \\ v_2
        \end{bmatrix}
    \end{equation}
    where $r > 0$ and $A_{\rm D}$ is the ideal diode defined by $A_{\rm D}(v) = \{0\}$ if $v<0$; $A_{\rm D}(0) = [0, + \infty)$; and $A_{\rm D}(v)=\emptyset$ otherwise. Define $\tilde{A}_{{\rm NPN},r} \coloneq A_{{\rm NPN},r} - [i_1\; i_2]^\top$. Then, any sequence $(v^k)_{k\in\mathbb{N}}$ satisfying the update rule \[
        v^{k+1} = T_{\gamma \tilde{A}_{{\rm NPN},r}}^B(v^k)
    \]
    with $B$ as in \Cref{cor:trans} and step size $\gamma >0$ converges weakly to $Bv^*$, where $v^*$ is a solution of \eqref{eq:transistor_eq}.
    
    Crucially, our iteration converges with any step size $\gamma>0$, whereas that of \cite[Prop.\,5.1]{quan2024scaled} only works for $\gamma > r(\sqrt{2} - 1)$, which is a restrictive condition, especially when $r\gg 1$, i.e., the case that usually occurs in practice. \Cref{fig:leaky_trans} shows the experimental results for the same configuration as in \cite[Fig.\,6]{quan2024scaled}.
\end{example}
\begin{proof}
    By \Cref{cor:trans}, we find that $(A_{\rm NPN}, B) \in \cM_0$, while per \Cref{def:semimon} with $\mu=\rho=0$, we also obtain that $((1/r)\id, B) \in \cM_0$. Similar to \Cref{ex:sum}, we can use \Cref{thm:scaling_prop:trans} to find $(A_{\rm NPN} + (1/r)\id, B) = (A_{{\rm NPN}, r}, B) \in \cM_0$. Since adding a constant does not change incremental properties, we have $(\tilde{A}_{{\rm NPN},r}, B)\in\cM_0$. From \Cref{prop:transf_resv}, we find that $T_{\gamma \tilde{A}_{{\rm NPN},r}}^B$ is firmly nonexpansive. It follows from \cite[Cor.\,5.17]{bauschke2017correction} that $(v^k)_{k\in\mathbb{N}}$ converges weakly to $\fix T_{\gamma \tilde{A}_{{\rm NPN},r}}^B = B(\zer \tilde{A}_{{\rm NPN},r})$.
\end{proof}

We next show how to leverage \Cref{prop:genprox} to derive a linearly preconditioned proximal point algorithm based on the transformed resolvent, before applying this to a (nonmonotone) common-emitter amplifier circuit.

Let $A:\cH\rightrightarrows\cH$ and suppose $\sing : \cH \to \cH$ is single-valued with $(A,\sing)\in\cM_0$. Let $M$ be a bounded, invertible linear operator on $\cH$. \Cref{prop:genprox} ensures that $((M^{-1})^*AM^{-1}, M\sing M^{-1})\in\cM_0$. We can then invoke \Cref{prop:transf_resv} and \KM\ \cite[Cor.\,5.17]{bauschke2017correction} to show that the iteration
\[
    x^{k+1} = M\sing M^{-1} \circ ((M^{-1})^* A M^{-1} + M\sing M^{-1})^{-1} x^k
\]
converges weakly to a point in $M\sing M^{-1}(\zer M^{-\top} A M^{-1})$, provided it exists. Now, define $\bar x^k $ by $M\bar x^k = x^k$. With some algebra, the iteration is then equivalent to 
\[
    \bar x^{k+1} \in \sing \circ (A + M^* M \sing)^{-1} M^* M \bar x^k.
\]    
Conversely, starting from a positive definite $P$, by letting $M=P^{1/2}$ be the positive definite square root of $P$, we find that the iteration
\begin{equation} \label{eq:tpppa}
     \bar x^{k+1} \in \sing \circ (A + P \sing)^{-1} P \bar x^k
\end{equation}
converges weakly to a point in $\sing(\zer A)$ if one exists. Further, instead of Fejér monotonicity in the Euclidean norm \cite[Cor.\,5.17(i)]{bauschke2017correction}, one now obtains Fejér monotonicity in the matrix $P$-norm \cite[Eq.\,(5.2.6)]{horn2012matrix}.

We now apply this iteration to a common-emitter amplifier circuit shown in \Cref{fig:common-emitter-amplifier}.
\begin{figure}[t]
\centering \vspace{0.15cm}
    \includegraphics[width=0.48\linewidth]{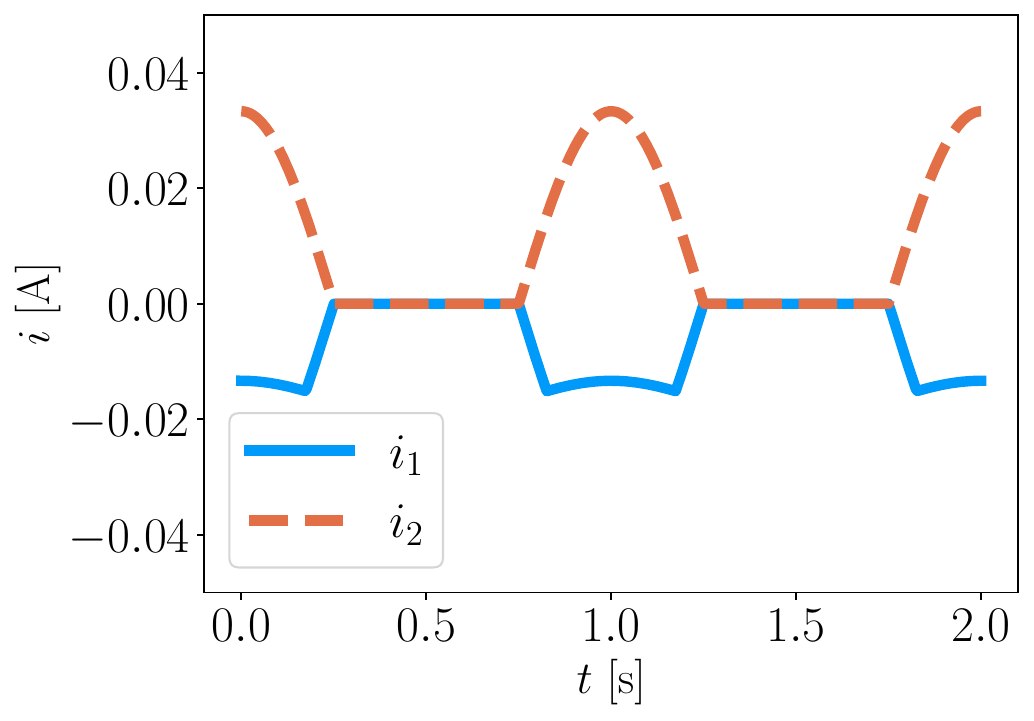}
    \hfill
    \includegraphics[width=0.48\linewidth]{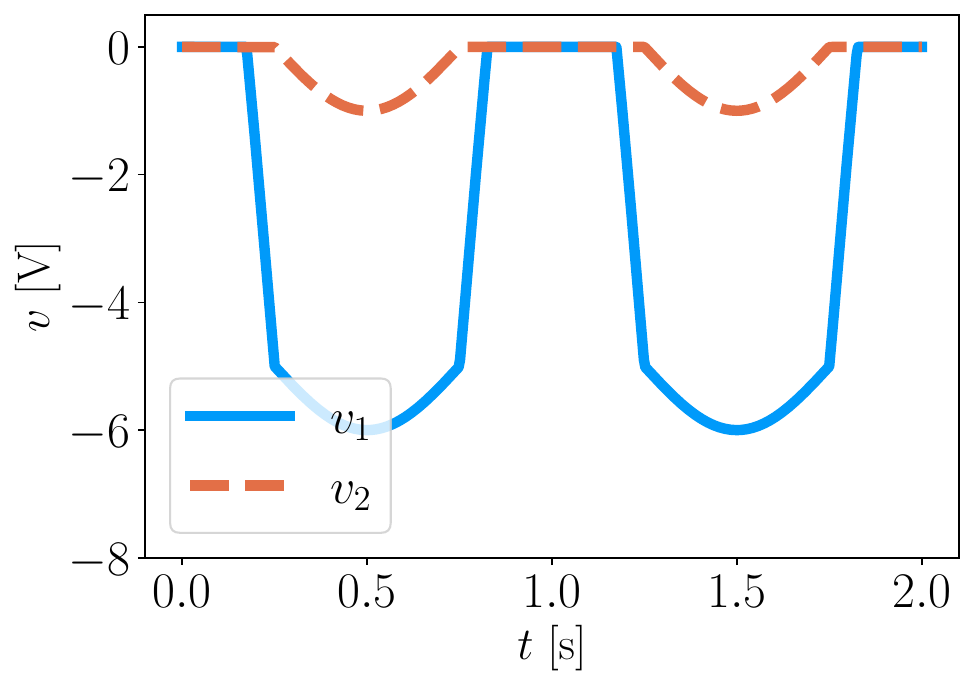}
    \caption{Solution to the inclusion problem \eqref{eq:cp_struc} for the common-emitter amplifier circuit in \cref{fig:common-emitter-amplifier} with linear resistors \(R_E = \SI{30}{\ohm}, R_C = \SI{300}{\ohm}\), and circuit parameters \(v_{\rm in} = \cos(2 \pi t) \si{\volt}\), \(v_+ = \SI{5}{\volt}\). These plots were obtained through the iteration \eqref{eq:chambolle-pock} with parameters \(\gamma = 10^{-3}, \tau = 100\).} \label{fig:amplifier}
\end{figure}
\begin{example} \label{ex:circuit}
    From \cite[Eq.\,(9)]{quan2024scaled}, the behavior of a common-emitter amplifier circuit can be obtained by solving an inclusion problem of the form
    \begin{equation} \label{eq:cp_struc}
        0 \in \begin{bmatrix}
            A(i) \\ B(v)
        \end{bmatrix} + \begin{bmatrix}
            0 & I_2^\top \\ -I_2 & 0
        \end{bmatrix}\begin{bmatrix}
            i \\ v
        \end{bmatrix} + \begin{bmatrix}s_v \\ s_i\end{bmatrix}
    \end{equation}
    where
    \[
        A \coloneq \begin{bmatrix}
            R_C & 0 \\ 0 & R_E
        \end{bmatrix}, \, B  \coloneq A_{\rm NPN}, \, s_v  \coloneq \begin{bmatrix}
            v_+ - v_{\rm in} \\ -v_{\rm in}
        \end{bmatrix}, \, s_i \coloneq 0
    \]
    with $v_{\rm in}\in\bR, v_+, R_C,R_E>0$ and $I_2$ denoting the identity matrix. Such a structure for solving circuits was studied in~\cite{chaffey2023circuit} in the monotone setting and in \cite{quan2024scaled} in the semimonotone setting. We now consider the pair monotonicity setting, that considers parameters that can not be covered by semimonotonicity.

    Let $R$ be associated with $A_{\rm NPN}$ as in \Cref{ex:srg-transistor} and suppose that $A^\top R^{-1} + R^{-\top}A \succeq 0$, $\tau, \gamma > 0$ and $\gamma\tau \|R\|^2 < 1$, where $\|R\|$ denotes the spectral norm of $R$. Consider the iteration
    \begin{equation} \label{eq:chambolle-pock}
        \begin{aligned}
            \bar i^k &= (A + \gamma^{-1}R^{-1})^{-1} (\gamma^{-1}i^k - R^\top v^k - s_v)\\
            \bar v^k &= (B + \tau^{-1}R^{-\top})^{-1} ( - Ri^k + \tau^{-1}v^k + 2\bar i^k - s_i) \\
           &\quad i^{k+1} = R^{-1} \bar i^k, \qquad v^{k+1} = R^{-\top} \bar v^k.
        \end{aligned}
    \end{equation}
    Then, $(Ri^k, R^\top v^k)_{k\in\mathbb{N}}$ converges weakly to a solution of \eqref{eq:cp_struc}, provided a solution exists.

    \Cref{fig:amplifier} shows the result of applying iteration \eqref{eq:chambolle-pock} to the common-emitter amplifier circuit in \Cref{fig:common-emitter-amplifier}. 
    For these numerical values, the required conditions hold. We emphasize that this setting is not covered by the Chambolle--Pock method \cite[Prop.\,5.2]{quan2024scaled}.
\end{example}
\begin{proof} 
    First, note that $(B, R^{-\top}), (A, R^{-1}) \in \cM_0$ where the first follows by  \Cref{cor:trans} and the second by the assumption and \Cref{def:semimon}. Similarly, the skew-symmetric term in \eqref{eq:cp_struc} can also be shown to be monotone with respect to $R^{-1} \oplus R^{-\top}$.
    
    By using the sum rule, \Cref{thm:scaling_prop:trans}, and the fact that constant terms do not affect the incremental properties, it follows that the complete operator in \eqref{eq:cp_struc} is monotone with respect to $\sing \coloneq (R^{-1} \oplus R^{-\top})$. Therefore, we can apply \eqref{eq:tpppa} with a positive definite $P$. Similar to how Chambolle--Pock \cite{chambolle2011first} arises from a specific choice of $P$ in the classical preconditioned proximal point algorithm to decouple the equations \cite[Eq.\,(1.1)]{evens2025convergenceCP}, we propose a similar preconditioner $P$, noting that our setting is different due to $\sing$.
    
    Let $P=\begin{bsmallmatrix}
         \gamma^{-1}I_2 & -R^\top \\ -R & \tau^{-1}I_2
    \end{bsmallmatrix}$, which is positive definite by the Schur complement \cite[Thm.\,7.7.7]{horn2012matrix} and our assumptions on $\tau,\gamma$. The iteration \eqref{eq:tpppa} applied to \eqref{eq:cp_struc} then becomes \eqref{eq:chambolle-pock} after some algebraic manipulations. The convergence then follows from the discussion above.
\end{proof}

\section{Conclusion}
In this paper, we introduced a novel scaled relative graph for pairs of operators. This framework naturally provides the geometric counterpart for recently introduced assumptions of paired monotonicity of operators. We have shown the practical relevance of these properties by computing the response of two highly nonmonotone circuits, thus extending known theoretical guarantees.

We believe that the proposed scaled relative graph for pairs of operators may prove valuable for stability analysis of feedback systems, and may simplify further analysis of other classes of nonmonotone circuits that can be handled by tailored splitting methods. Other interesting avenues for future work include developing numerical methods for calculating paired SRGs and nonlinear resolvents.

\begin{appendix}
For the proofs, we require the following equivalent formulation of the complex-conjugate pair $z_{\pm}$ (see e.g., \cite[Eq.\,(1)]{ryu2022scaled}): 
\begin{equation} \label{eq:proj}
    z_{\pm}(u - \bar u, x - \bar x) = \frac{\langle u-\bar u,x-\bar x \rangle}{\|x-\bar x\|^2} \pm i\frac{\|\proj_{\{x-\bar x\}^\perp}(u-\bar u)\|}{\|x-\bar x\|},
\end{equation}
where $\proj_{\{x-\bar x\}^\perp}$ is the projection onto the subspace orthogonal to $x-\bar x$.

\noindent \textbf{Proof of \cref{thm:scaling_prop}}
\begin{proof}
    ``\labelcref{thm:equiv}'': Denote $C\coloneq A\circ B^{-1}$, which has graph given by $\graph C=\{(u_B,u_A)\in\cH\times\cH:\exists x\in\cH,\ (x,u_B)\in \graph B,\ (x, u_A)\in \graph A\}$. Let $z\in\cG(A,B)\setminus\{\infty\}$. Then there exist $(x,u_A),(\bar x, \bar u_A)\in \graph A$ and $(x, u_B),(\bar x, \bar u_B)\in \graph B$ with $u_B\neq\bar u_B$, hence $z=z_\pm(u_A-\bar u_A,u_B-\bar u_B)$. Since $(u_B,u_B),(\bar u_B,\bar u_B)\in\graph\id$ and $(u_B,u_A),(\bar u_B,\bar u_A)\in\graph C$, the same point $z$ belongs to $\cG(C,\id)$.
    
    Conversely, if $z\in\cG(C,\id)\setminus\{\infty\}$, then $z=z_\pm(u_A-\bar u_A,u_B-\bar u_B)$ for some $(u_B,u_A),(\bar u_B,\bar u_A)\in\graph C$ with $u_B\ne\bar u_B$. Therefore, there exist $x,\bar x\in\cH$ such that $(x,u_B), (\bar x, \bar u_B)\in \graph B$, $(x,u_A), (\bar x, \bar u_A)\in \graph A$. Hence \(z\in\cG(A,B)\), and the finite parts coincide.
    
    It remains to check the point at infinity. By \Cref{def:pair_srg}, $\infty\in\cG(A,B)$ if and only if there exist $(x,u_A), (\bar x, \bar u_A)\in \graph A$ and $(x,u_B),(\bar x, \bar u_B)\in \graph B$ such that $u_B=\bar u_B$ and $u_A\neq\bar u_A$. Denote $y=u_B=\bar u_B$ and note that $(y,u_A),(y,\bar u_A)\in\graph C$, with distinct outputs at the same input $y$, implying $\infty\in\cG(C,\id)$. Conversely, if $\infty\in\cG(C,\id)$, then there exist $(y,u_A),(y,\bar u_A)\in\graph C$ with $u_A\neq\bar u_A$. Thus, there exist $x,\bar x\in\cH$ such that $(x,y),(\bar x, y)\in \graph B$, $(x,u_A),(\bar x, \bar u_A)\in \graph A$. Taking $u_B=\bar u_B=y$ gives $\infty\in\cG(A,B)$, and we conclude that $\cG(A,B)=\cG(A\circ B^{-1},\id)$.

    ``\labelcref{thm:scaling_prop:scale}'': This follows from \cref{def:pair_srg} along with the properties of the inner product and the norm.

    ``\labelcref{thm:scaling_prop:inv}'': The equality follows for points that are neither $0$ nor $\infty$ by definition. If $\infty \in \cG(A,B)$, there exist $(x,u_A),(\bar x, \bar u_A) \in \graph A$ and $(x,u_B),(\bar x,\bar u_B) \in \graph B$ such that $u_B = \bar u_B$ and $u_A \neq \bar u_A$. This immediately implies that $0 \in \cG(B,A)$. The zero case follows similarly.

    ``\labelcref{thm:scaling_prop:id}'': The first equality is \Cref{thm:scaling_prop:inv} and the second is \cite[Thm.\,5]{ryu2022scaled} since $\cG(A, \id)$ is the standard SRG.

    ``\labelcref{thm:scaling_prop:trans}'': The inclusion for the finite points follows similarly to \cite[Thm.\,6]{ryu2022scaled}. If $\infty\in\cG(A+B,C)$, then there exist $(x,u_C),(\bar x,\bar u_C)\in\graph C$, $(x,u_A),(\bar x,\bar u_A)\in\graph A$, and $(x,u_B),(\bar x,\bar u_B)\in\graph B$ such that $u_C=\bar u_C$ and $u_A+u_B\neq\bar u_A+\bar u_B$. Hence $u_A\ne\bar u_A$ or $u_B\ne\bar u_B$. Therefore $\infty\in\cG(A,C)$ or $\infty\in\cG(B,C)$ and hence $\infty \in \cG(A,C) + \cG(B,C)$.

    ``\labelcref{thm:scaling_prop:trans_inv}'': Combine \Cref{thm:scaling_prop:trans} with \Cref{thm:scaling_prop:inv}, and use that the inverse is defined elementwise, so subset inclusions are preserved.
 \end{proof}

\noindent\textbf{Proof of \cref{thm:single_valued}}
\begin{proof}
    ``\labelcref{thm:single_valued:id}'': Take $(x, \sing(x)), (\bar x, \sing(\bar x)) \in \graph \sing$. Clearly, $\infty \notin \cG(\sing,\sing)$ since $\sing$ is a function and we would require both $F(x) = F(\bar x)$ and $F(x)\neq F(\bar x)$. Since $F$ is not constant, there exist $x,\bar x$ with $F(x) \neq F(\bar x)$. For any such pair, we have that $z_{\pm} = \tfrac{\|\sing(x)-\sing(\bar x)\|}{\|\sing(x)-\sing(\bar x)\|} \exp(\pm i0) = \{1\}$. Thus $\cG(\sing,\sing) = \{1\}$.
    
    ``\labelcref{thm:single_valued:sum}'': Take $z \in \cG(A+\sing,\sing) \setminus \{\infty\}$. Then, there exist $x, \bar x$ and $(x, u_A),(\bar x, \bar u_A) \in \graph A$ such that $z \in z_{\pm}(\sing(x)+u_A-(\sing(\bar x) + \bar u_A), \sing(x)-\sing(\bar x))$. Since $z \neq \infty$, $\sing(x) \neq \sing(\bar x)$ and $z_{\pm}(\sing(x)+u_A-(\sing(\bar x) + \bar u_A),\sing(x)-\sing(\bar x)) = 1 + \frac{\langle u_A-\bar u_A,\sing(x)-\sing(\bar x) \rangle}{\|\sing(x)-\sing(\bar x)\|^2} \pm i\frac{\|\proj_{\{\sing(x)-\sing(\bar x)\}^\perp}(u_A-\bar u_A)\|}{\|\sing(x)-\sing(\bar x)\|} \subseteq 1+\cG(A,\sing)$. Now, $\infty \in \cG(A+\sing,\sing)$ implies $\sing(x) = \sing(\bar x)$ and $u_A + \sing(x) \neq \bar u_A + \sing(\bar x)$, i.e. $u_A \neq \bar u_A$ which then means that $\infty \in \cG(A,\sing)$. The opposite inclusions follow similarly.

    ``\labelcref{thm:single_valued:sum_inv}'': Apply \Cref{thm:equiv} to the pair $(\sing, A + \sing)$.

    ``\labelcref{thm:single_valued:precomp}'': Take $z \in \cG(A \circ \sing, B \circ \sing) \setminus \{\infty\}$. Then, there exist $(x, u), (\bar x, \bar u) \in \graph (A \circ \sing)$ and $(x,v),(\bar x, \bar v) \in \graph (B \circ \sing)$ such that $z \in z_{\pm}(u-\bar u,v-\bar v )$ and $v \neq \bar v$. This implies $(\sing(x),u),(\sing(\bar x),\bar u) \in \graph A$ and $(\sing(x), v),(\sing(\bar x), \bar v) \in \graph B$ with $v \neq \bar v$. Then, by definition, $z \in \cG(A,B)$. Now take $z = \infty$. Then, for some pair as before, $v = \bar v$ and $u \neq \bar u$ meaning also that $\infty \in \cG(A,B)$.

    Now assume, moreover, that $\sing$ is surjective. Take $z \in \cG(A,B) \setminus \{\infty\}$. Then, there exist $(x, u_A), (\bar x, \bar u_A) \in \graph A$ and $(x,u_B),(\bar x, \bar u_B) \in \graph B$ such that $z \in z_{\pm}(u_A-\bar u_A,u_B-\bar u_B)$ and $u_B \neq \bar u_B$. Clearly, since $\sing$ is surjective, there exist $y, \bar y \in \cH$ such that $\sing(y) = x$ and $\sing(\bar y) = \bar x$. Then, $(y, u_A),(\bar y, \bar u_A) \in \graph (A \circ F)$ and $(y, u_B),(\bar y, \bar u_B) \in \graph (B \circ F)$ with $u_B \neq \bar u_B $ and thus $z \in z_{\pm}(u_A-\bar u_A, u_B-\bar u_B) \subset \cG(A \circ \sing, B \circ \sing)$. The case $z = \infty$ follows similarly.

    ``\labelcref{thm:single_valued:lip}'': Take $z \in \cG(\sing \circ A, B)$. If $z = \infty$, then there exist $(x, u),(\bar x, \bar u) \in \graph (\sing \circ A)$ and $(x,v),(\bar x, \bar v ) \in \graph B$ such that $v = \bar v$ and $u \neq \bar u$. Therefore, there also exist $(x,y),(\bar x,\bar y) \in \graph A$ such that $u = \sing(y), \bar u = \sing(\bar y)$ and, since $\sing$ is single-valued and $u \neq \bar u$, we have that $y \neq \bar y$. This implies $\infty \in \cG(A,B)$, which would contradict the assumption that $\cG(A, B) \subseteq D(0,l)$, and therefore, $\infty\notin \cG(\sing \circ A, B)$. Now assume that $z \neq \infty$. This means that with pairs as before such that $v \neq \bar v$ we have $|z| = \tfrac{\|u-\bar u\|}{\|v-\bar v\|} = \frac{\|\sing(y)-\sing(\bar y)\|}{\|v -\bar v\|} \leq L\frac{\|y-\bar y\|}{\|v -\bar v\|}$, since $\sing$ is a $L$-Lipschitz operator in light of \cite[Prop.\,1 and Thm.\,2]{ryu2022scaled}. This implies that $|z| \leq L|\bar z|$ for $\bar z \in \cG(A,B)$ and thus that $|z| \leq L\sup\{|\bar z|:\bar z \in \cG(A,B)\}$. Lastly, from the hypothesis, $\sup\{|\bar z|:\bar z \in \cG(A,B)\} \leq l$, so we obtain the claimed result $\cG(\sing\circ A, B) \subseteq D(0,Ll)$.

    ``\labelcref{thm:single_valued:lip2}'': Take $z\in\cG(A, \id)$. If $z=\infty$, then there exist $(x,u),(\bar x,\bar u) \in \graph A$ such that $u\neq \bar u$ and $x=\bar x$. Since $\sing$ is single-valued, this implies that there exist $(x,u),(\bar x,\bar u) \in \graph A$ such that $u \neq \bar u$ and $\sing(x) = \sing(\bar x)$, i.e., $\infty \in \cG(A, \sing)$, which contradicts the assumption that $\cG(A,\sing)\subseteq D(0,L)$, so $\infty\notin\cG(A,\id)$. In particular, we have that $A$ is single-valued.
    
    Further, for $z\neq\infty$, as in the proof of \Cref{thm:single_valued:lip}, we know that $\sing$ is an $l$-Lipschitz operator. Similarly, from $\cG(A,\sing) \subseteq D(0,L)$ we have that $z' \in \cG(A,\sing)$ implies $|z'| \leq L$ and thus that $\|A(x)- A(\bar x)\| \leq L\|\sing(x) - \sing(\bar x)\|$ for all $x,\bar x\in\dom A$. Then, if $z \in \cG(A, \id)$ we have that $|z|=\frac{\|A(x)- A(\bar x)\|}{\|x-\bar x\|}$ for some $x,\bar x\in\dom A$. Using the previous inequality along with the Lipschitz continuity of $\sing$ we obtain $|z| \leq lL$, which is the claimed result.

    ``\labelcref{thm:single_valued:linmap}'': To begin with, note that since $M$ is bounded and invertible, we have that $(M^*)^{-1} = (M^{-1})^*$. For any $z\in\cG(A,B), z \in \bC_{\geq 0}$. Thus, for all $(x,u_A),(\bar x,\bar u_A) \in \graph A$ and $(x,u_B),(\bar x,\bar u_B) \in \graph B$ such that $u_B\neq\bar u_B$ we have that $\langle u_A - \bar u_A, u_B - \bar u_B\rangle \geq 0$ or, equivalently, $\langle (M^*)^{-1}(u_A - \bar u_A), M(u_B - \bar u_B)\rangle \geq 0$ by definition of the adjoint \cite[Thm.\,4.10]{Rudin1991}. Since $M$ is invertible, it holds for any $(x,u),(\bar x,\bar u)\in\graph((M^*)^{-1} \circ A)$ and $(x,v),(\bar x,\bar v)\in\graph (M\circ B)$ with $v\neq\bar v$ that $\langle u-\bar u, v - \bar v\rangle \geq 0$. The $\infty$ case follows similarly and we conclude that $\cG((M^{-1})^* \circ A, M \circ B) \subseteq \bC_{\geq0}$.
\end{proof}

\noindent\textbf{Proof of \Cref{prop:srgfull}}
\begin{proof}
    Suppose that there exists a nonnegatively homogeneous function $h$ that certifies membership as in \Cref{prop:srgfull}. Notice that $(A, B) \in \cP \implies \cG(A, B) \subseteq \cG(\cP)$ by construction, and it remains to show that $\cG(A, B) \subseteq \cG(\cP) \implies (A, B) \in \cP$. For all points $(x, u_A), (\bar x, \bar u_A) \in\graph A$ and $(x,u_B), (\bar x, \bar u_B)\in\graph B$ such that $u_B\neq\bar u_B$, the same technique from \cite[Thm.\,2]{ryu2022scaled} yields that $\cG(A, B) \subseteq \cG(\cP) \implies h(\|u_A - \bar u_A\|^2, \|u_B - \bar u_B\|^2, \langle u_A - \bar u_A, u_B - \bar u_B\rangle) \leq 0$.

    Now suppose $u_B = \bar u_B$, while $u_A \neq \bar u_A$. Then $\infty \in \cG(A, B) \subseteq \cG(\cP)$. Thus, there is some $(C, D) \in \cP$ such that $(x', u_C'), (\bar x', \bar u_C') \in \graph C$, $(x', u_D'), (\bar x', \bar u_D') \in \graph D$, and $u_D' = \bar u_D'$, $u_C' \neq \bar u_C'$. Since $(C, D) \in \cP$, we have that $h(\|u_C' - \bar u_C'\|^2, 0, 0) \leq 0$. Multiplying by $\|u_A - \bar u_A\|^2/\|u_C' - \bar u_C'\|^2$ and using homogeneity, we find $h(\|u_A - \bar u_A\|^2, 0,0) \leq 0$. Lastly, if $u_B =\bar u_B$ and $u_A=\bar u_A$, then $h$ is evaluated at all zeros, and its output must therefore be zero by homogeneity. In all three cases, we have $h(\|u_A - \bar u_A\|^2, \|u_B - \bar u_B\|^2, \langle u_A - \bar u_A, u_B - \bar u_B\rangle) \leq 0$ and since the points were arbitrary evaluations of $A,B$, we have that $(A,B)\in \cP$ by the hypothesis.

    Conversely, suppose that $\cP$ is SRG-full. Let $K = \{(a, b, c) \in \bR^3 \mid a\geq 0, b\geq 0, c^2 \leq ab\}$ and consider
    \[
        h(a,b,c) = \begin{cases}
            (a+b) I(\tilde{z}(a, b, c)) & \textnormal{if $(a,b,c) \in K \setminus \{0\}$}, \\
            0 & \textnormal{otherwise}
        \end{cases}
    \]
    where $\tilde{z}(a,b,c) = \frac{c}{b} + i\frac{\sqrt{ab-c^2}}{b}$ if $b > 0$ and $\tilde{z}(a,b,c) = \infty$ if $b=0$ and $a > 0$, and
    \[
        I(z) = \begin{cases}
            0 & \textnormal{if $z\in\cG(\cP)$}, \\
            1 & \textnormal{otherwise}.
        \end{cases}
    \]

    By our constraint on $(a,b,c)$, $\tilde{z}$ will not be evaluated on points where it is not defined. Observe also that $\tilde{z}(\|u - \bar u\|^2, \|x - \bar x\|^2, \langle u-\bar u, x - \bar x\rangle) = z_+(u - \bar u, x - \bar x)$ as in \eqref{eq:proj}.

    Let $\kappa > 0$ and $(a,b,c) \in K \setminus\{0\}$. Then also $(\kappa a, \kappa b, \kappa c) \in K \setminus \{0\}$ and $\tilde{z}(\kappa a, \kappa b, \kappa c) = \tilde{z}(a, b, c)$ by definition. It follows that $h(\kappa a, \kappa b, \kappa c) = (\kappa a+\kappa b) I(\tilde{z}(a,b,c)) = \kappa h(a,b,c)$. The cases $\kappa = 0$ and $(a,b,c) \notin K \setminus\{0\}$ both yield $h(\kappa a, \kappa b, \kappa c) = 0 = \kappa h(a,b,c)$, so we obtain that $h$ is nonnegatively homogeneous. It remains to show that $(A,B) \in \cP \iff h(\|u_A - \bar u_A\|^2, \|u_B - \bar u_B\|^2, \langle u_A - \bar u_A, u_B - \bar u_B\rangle) \leq 0$ for all $(x, u_A), (\bar x, \bar u_A)\in\graph A$ and $(x,u_B), (\bar x, \bar u_B) \in \graph B$.

    Thus, first assume that $(A,B) \in \cP$. Let $(x,u_A), (\bar x, \bar u_A) \in \graph A$ and $(x,u_B), (\bar x, \bar u_B)\in \graph B$ and denote $a = \|u_A - \bar u_A\|^2$, $b = \|u_B - \bar u_B\|^2$ and $c = \langle u_A - \bar u_A, u_B - \bar u_B\rangle$. Note that $(a,b,c) \in K$ by Cauchy--Schwarz. If $(a,b,c) = 0$, then $h(a,b,c) = 0$ by construction. Otherwise, if $(a,b,c)\neq 0$, \Cref{def:pair_srg} and \eqref{eq:proj} ensure that $\tilde{z}(a,b,c) \in \cG(A,B)$. Since $\cP$ is SRG-full, it follows that $\tilde{z}(a,b,c) \in \cG(\cP)$ and $h(a,b,c) = (a+b) \cdot 0 = 0$ by construction. Thus, $h(a,b,c) \leq 0$ for all points in the graphs.

    Second, assume $h(\|u_A - \bar u_A\|^2, \|u_B - \bar u_B\|^2, \langle u_A - \bar u_A, u_B - \bar u_B\rangle) \leq 0$ for all $(x, u_A), (\bar x, \bar u_A)\in\graph A$ and $(x,u_B), (\bar x, \bar u_B) \in \graph B$. It suffices to show that $\cG(A,B)\subseteq \cG(\cP)$ from which the desired result follows by SRG-fullness of $\cP$. To this end, let $z\in\cG(A,B)$. Assume $z\neq\infty$. Since the scaled relative graph is symmetric with respect to the real axis, we may assume that $\impart z \geq 0$. By \Cref{def:pair_srg}, there exists some $(x, u_A), (\bar x, \bar u_A)\in\graph A$ and $(x,u_B), (\bar x, \bar u_B) \in \graph B$ so that $z = \tilde{z}(a, b, c)$ and $a + b > 0$, where we denoted $a = \|u_A - \bar u_A\|^2$, $b = \|u_B - \bar u_B\|^2$ and $c = \langle u_A - \bar u_A, u_B - \bar u_B\rangle$. Further, since $h$ is nonnegative by construction, $h(a,b,c)\leq 0 \iff h(a,b,c) = 0$. Therefore, since $a+b > 0$, we require $I(\tilde{z}(a,b,c)) = 0$, implying that $z = \tilde{z}(a,b,c) \in \mathcal{G}(\cP)$ by definition of $I$. For $z=\infty$, use $a > 0, b=0, c=0$, so that $\tilde{z}(a,b,c) =\infty$. Since $h(a,b,c)\leq 0$, $h\geq 0$, and $a+b>0$, we obtain $I(\infty)=0$, hence $\infty\in\cG(\cP)$.
\end{proof}

\noindent\textbf{Proof of \Cref{prop:srg_semi}}
\begin{proof}
    That $\cS_{\mu,\rho}$ is SRG-full follows from \Cref{def:semimon} and \Cref{prop:srgfull} with $h : (a,b,c)\mapsto \rho a + \mu b - c$.

    We now show the claimed expression for $\cG(\cS_{\mu,\rho})$. Suppose that $(A, B) \in \cS_{\mu,\rho}$, $(x, u_A), (\bar x, \bar u_A) \in \graph A$ and $(x, u_B), (\bar x, \bar u_B) \in \graph B$ are arbitrary with $u_B\neq\bar u_B$. It follows from \Cref{def:semimon} that $\langle u_A - \bar u_A, u_B - \bar u_B\rangle \geq \mu \|u_B - \bar u_B\|^2 + \rho \|u_A - \bar u_A\|^2$, and dividing by $\|u_B - \bar u_B\|^2$ yields $\frac{\langle u_A - \bar u_A, u_B - \bar u_B\rangle}{\|u_B - \bar u_B\|^2} \geq \mu + \rho \frac{\|u_A - \bar u_A\|^2}{\|u_B - \bar u_B\|^2}$ which is equivalent to  $\realpart(z_{\pm}(u_A - \bar u_A, u_B - \bar u_B)) \geq \mu + \rho |z_{\pm}(u_A - \bar u_A, u_B - \bar u_B)|^2$, see \eqref{eq:proj}. Now suppose $\rho > 0$. Clearly, if $u_B = \bar u_B $ then also $u_A = \bar u_A$ from \Cref{def:semimon} and thus $\infty \notin \cG(\cS_{\mu,\rho})$.

    For the opposite inclusion, suppose $\rho=0$. Then the set on the right-hand side is exactly $\bC_{\geq\mu}$, and from \cite[Prop.\,1]{ryu2022scaled}, we find $\bC_{\geq\mu} = \cup_{(A, \id) \in \cM_\mu}\cG(A, \id) \subseteq \cG(\cM_\mu) = \cG(\cS_{\mu, 0})$. Otherwise, if $\rho> 0$, let $z  = x + iy$ satisfy $\realpart(z) \geq \mu + \rho|z|^2$. Then $\rho x^2 - x + \rho y^2 + \mu \leq 0$. Further, dividing by $\rho$ and completing the square, we obtain
    \[
        \left(x - \frac{1}{2\rho}\right)^2 + y^2 \leq\left(\frac{1}{2\rho}\right)^2 - \frac{\mu}{\rho} = \frac{1-4\mu\rho}{4\rho^2}.
    \]
    Observe that this defines the same disk as in \cite[Prop.\,3.4]{quan2024scaled}. Therefore, $\cup_{(A, \id) \in \cS_{\mu,\rho}} \cG(A, \id) \subseteq \cG(\cS_{\mu,\rho})$. The case $\rho < 0$ follows similarly.
\end{proof}
\end{appendix}

\newpage
\bibliographystyle{plain}
\bibliography{references}

\end{document}